\documentclass[preprint,11pt]{elsarticle}

\usepackage{amssymb}
\usepackage{amsmath}
\usepackage{amsfonts}
\usepackage{amssymb}
\usepackage{indentfirst,latexsym,bm}
\usepackage{amsthm}
\usepackage{color,xcolor}
\usepackage[all]{xy}
\input{mathrsfs.sty}
\newtheorem{theorem}{Theorem}[section]
\newtheorem{lemma}[theorem]{Lemma}
\newtheorem{corollary}[theorem]{Corollary}

\newtheorem{problem}[theorem]{Problem}

\theoremstyle{definition}
\newtheorem{definition}{Definition}[section]

\theoremstyle{definition}

\theoremstyle{remark}
\newtheorem{remark}{Remark}[section]
\theoremstyle{question}

\numberwithin{equation}{section}

\journal{XXX}

\begin{document}

\begin{frontmatter}





\title{The Frobenious distances from projections to an idempotent matrix}
\author[shnu]{Xiaoyi Tian}
\ead{tianxytian@163.com}
\author[shnu]{Qingxiang Xu\corref{cor2}}
\ead{qingxiang\_xu@126.com}
\author[HS]{Chunhong Fu}
\ead{fchlixue@163.com}
\address[shnu]{Department of Mathematics, Shanghai Normal University, Shanghai 200234, PR China}
\address[HS]{Health School Attached to Shanghai University of Medicine $\&$ Health Sciences,   Shanghai 200237, PR China}

\begin{abstract} For each pair of matrices $A$ and $B$ with the same order, let
$\|A-B\|_F$ denote their Frobenius distance. This paper deals mainly with the Frobenius  distances from projections to an idempotent matrix.
For every idempotent $Q\in \mathbb{C}^{n\times n}$, a projection $m(Q)$ called the matched projection can be induced. It is proved that $m(Q)$ is the unique projection whose Frobenius  distance away from $Q$ takes the minimum value among all the Frobenius  distances from projections to $Q$, while $I_n-m(Q)$ is the unique projection whose Frobenius distance away from $Q$ takes the maximum value. Furthermore, it is proved that for every number $\alpha$ between the minimum value and the maximum value, there exists a projection $P$ whose Frobenius  distance away from $Q$ takes the value $\alpha$. Based on the above characterization of the minimum distance, some Frobenius norm upper bounds and lower bounds of $\|P-Q\|_F$ are derived under the condition of $PQ=Q$ on a projection $P$ and an idempotent $Q$.
  \end{abstract}

\begin{keyword}Idempotent, projection, Frobenius norm, norm bound
\MSC 65F35, 15A60, 15A09



\end{keyword}

\end{frontmatter}



\section{Introduction}\label{sec:pre}
By a projection, we mean a Hermitian idempotent.  The direct sum decompositions and the orthogonal decompositions of spaces are frequently used, so
idempotents and projections (also known as skew projections and orthogonal projections) have been intensely studied, and various applications can be found in the literature. Yet, there are still some fundamental issues that remain to be unknown. One of which is concerned with the Frobenius distances from projections to an idempotent, which will be addressed in the following Problem~\ref{problems}.

Let   $\mathbb{C}^{m\times n}$ be the set of   $m\times n$ complex matrices, and let $\mathcal{P}(\mathbb{C}^{n\times n})$ be the set consisting of all projections in $\mathbb{C}^{n\times n}$. The identity matrix and the zero matrix in $\mathbb{C}^{n\times n}$ are denoted by $I_n$ and $0_n$, respectively. Given a pair of matrices $A$ and $B$ in $\mathbb{C}^{m\times n}$, the Frobenius norm $\|A-B\|_F$ is referred to be the  (Frobenius) distance between $A$ and $B$. Suppose now that $Q\in \mathbb{C}^{n\times n}$ is an idempotent. Let
\begin{align}\label{defn of ran Q}&\mbox{ran}_Q=\big\{\|P-Q\|_F: P\in \mathcal{P}(\mathbb{C}^{n\times n})\big\},\\
\label{Q min}&\mbox{min}_Q=\inf\big\{\lambda:\lambda\in \mbox{ran}_Q\big\},\quad \mbox{max}_Q=\sup\big\{\lambda:\lambda\in \mbox{ran}_Q\big\}.
\end{align}
A couple of  issues associated with a general idempotent $Q$ can be raised naturally as follows.
\begin{problem}\label{problems}Suppose that $Q\in \mathbb{C}^{n\times n}$ is an idempotent such that $Q\ne I_n$ and $Q\ne 0_n$.
\begin{enumerate}
\item[{\rm (i)}] Is it possible to determine the numbers $\mbox{min}_Q$ and $\mbox{max}_Q$?
\item[{\rm (ii)}] Is it possible to find out projections $P_1$ and $P_2$ such that
$$\|P_1-Q\|_F=\mbox{min}_Q\quad\mbox{and}\quad \|P_2-Q\|_F=\mbox{max}_Q?$$
\item[{\rm (iii)}] If the answer to (ii) is positive, whether such projections $P_1$ and $P_2$ are unique?
\item[{\rm (iv)}] If the answer to (ii) is positive, whether it is true that $[\mbox{min}_Q,\mbox{max}_Q]=\mbox{ran}_Q$?
\end{enumerate}
\end{problem}

Evidently, the same problems are also worthwhile to be investigated in terms of  the operator norm (spectral norm). Actually, in the latter case some partial
answers can be found in \cite[Section~4]{TXF}. Specifically, a new term called the matched projection is introduced recently in \cite[Section~3]{TXF} for an arbitrary idempotent on a Hilbert $C^*$-module.

To get a deeper understanding of the distances from projections to an idempotents, in this paper we restrict our attention to the matrix case and adopt the Frobenius norm instead of the operator norm. Given an arbitrary idempotent $Q\in \mathbb{C}^{n\times n}$, let $m(Q)$ be its matched projection defined  as \eqref{equ:final exp of P}.
In Lemma~\ref{lem:f norm m Q minus Q} of this paper,  a formula for $\|m(Q)-Q\|_F$ is derived . It is proved in Theorem~\ref{thm:existence} that
 $\|m(Q)-Q\|_F\le \|P-Q\|_F$ for every projection $P\in \mathbb{C}^{n\times n}$, and it is shown furthermore in Theorem~\ref{thm:uniqueness} that the equality above occurs only if $P=m(Q)$.
Thus, the optimal approximation from projections to an idempotent is  clarified.

Based on the above characterization of  the optimal approximation, in Section~\ref{sec:additional remark} of this paper we will provide the positive answers to all the remaining issues that stated in Problem~\ref{problems}; see Theorems~\ref{thm:maximum distance} and \ref{intermediate value theorem} for the details. So, a complete answer to Problem~\ref{problems} is carried out.

 Let $P\in \mathbb{C}^{n\times n}$ be a projection and $Q\in \mathbb{C}^{n\times n}$ be  an idempotent.  Under the condition of  $PQ=Q$, one lower bound and one upper bound  of $\|P-Q\|_F$ are derived in Theorems~\ref{thm:app lower-01} and \ref{thm:app up-01}, respectively.
As a result, an upper bound of $\|P-Q^\dag\|_F$ is also derived in Theorem~\ref{thm:app up-01 point 5}.
Note that the above condition is trivially satisfied when $P$ is taken to be $P_{\mathcal{R}(Q)}$ (the projection from $\mathbb{C}^{n}$ onto the range of $Q$).
It is particularly interesting to make a detailed comparison between $\|m(Q)-Q\|_F$ and $\|P_{\mathcal{R}(Q)}-Q\|_F$, which is dealt with in  Theorem~\ref{thm:the range projection} (see also Remark~\ref{rem:optimal numbers} for a supplement). Observe that the condition of $PQ=Q$ is also trivially satisfied when
$P=AA^\dag$ and $Q=AX$, in which $A^\dag$ denotes the Moore-Penrose inverse of $A$, while $X$ is an arbitrary inner inverse or outer inverse of $A$.
Hence, a unified way  can be employed to deal with upper bounds and lower bounds of $\|AA^\dag-AX\|_F$; see Corollaries~\ref{thm:drazin lower bound} and \ref{thm:drazin upper bound} in the special case that $X$ is the Drazin inverse of $A$.

The paper is organized as follows. Some basic knowledge about trace, the weak majorization and the matched projection are included in Section~\ref{sec:pre}.
 Section~\ref{sec:positive answer} is focused on the study of the minimum distance from projections to an idempotent,  while Section~\ref{sec:additional remark} is devoted to clarifying the maximum distance and the range of the distances from projections to an idempotent. As applications, some Frobenius norm bounds are derived in the last section.

\section{Some preliminaries}\label{sec:pre}
Throughout the rest of this paper, $\mathbb{N}$ is the set of
positive integers,  $\mathbb{C}_r^{m\times n}$ is the subset of  $\mathbb{C}^{m\times n}$ consisting of matrices with rank $r$. For each $A\in \mathbb{C}^{m\times n}$, its (column) range, conjugate transpose, Moore-Penrose inverse \cite{XWG},  Frobenius norm and spectral norm (2-norm) are denoted by $\mathcal{R}(A)$, $A^*$, $A^\dag$, $\Vert A\Vert_F$ and $\|A\|_2$, respectively. When $A$ is a square matrix, ${\rm tr}(A)$ stands for the trace of $A$.
Let $|A|$ be the square root of $A^*A$. The notation $P_{\mathcal{R}(Q)}$ is used
for an idempotent $Q\in \mathbb{C}^{n\times n}$ to denote the projection from $\mathbb{C}^n$ onto $\mathcal{R}(Q)$.

We begin with an auxiliary lemma, which contains several elementary known results.
\begin{lemma}\label{elemental observations}Let $A,B,T=(t_{ij})\in\mathbb{C}^{n\times n}$,  and let $x$ and $y$ be elements in an inner-product space.
\begin{enumerate}
\item[{\rm (i)}] Suppose that $A$ and $B$ are Hermitian such that  $A-B$ is positive semi-definite. Then $\mbox{tr}(A)\ge \mbox{tr}(B)$, and $\mbox{tr}(A)=\mbox{tr}(B)$ only if $A=B$.
\item[{\rm (ii)}] If $\sum\limits_{i=1}^n |t_{ii}|^2=\sum\limits_{i=1}^n \sigma_i^2(T)$, then $T$ is diagonal, where $\sigma_i(T) \,(1\le i\le n)$ denote the singular values of $T$.
\item[{\rm (iii)}] If $\langle x,y\rangle=\|x\|\cdot \|y\|$, then $\|x\|\cdot y=\|y\|\cdot x$.
\end{enumerate}
\end{lemma}
\begin{proof} (i) Write $(A-B)^\frac12$ simply as $C$. Then
$$\mbox{tr}(A)-\mbox{tr}(B)=\mbox{tr}(A-B)=\|C\|_F^2.$$
Hence, the conclusion follows.

(ii) The conclusion follows immediately from the following equations
$$\sum_{i=1}^n\sum_{j=1}^n|t_{ij}|^2=\mbox{tr}(TT^*)=\sum_{i=1}^n \sigma^2_i(T).$$

(iii)  Put $z=\|x\|\cdot y-\|y\|\cdot x$. From the assumption it is easily seen that $\langle z,z\rangle=0$, so the desired conclusion follows.
\end{proof}

Next, we recall some characterizations of the weak majorization associated with convex functions \cite[Chapter 10]{Zhang}.
Given $a =(a_1,  \cdots,a_n)$,  $b=(b_1,  \cdots,b_n) \in  \mathbb{R}^n$, let $a^\downarrow =(a_1^{\downarrow}, \cdots, a_n^{\downarrow})$ and $b^\downarrow=(b_1^{\downarrow},  \cdots, b_n^{\downarrow})$ be
the decreasing rearrangements of $a$ and $b$, respectively. Recall that $a$ is said to be weakly majorized by $b$, denoted by $a \prec_w b$, if
$$\sum_{i=1}^k a_i^\downarrow \leq \sum_{i=1}^k b_i^\downarrow \quad (1\leq k \leq n).$$

\begin{lemma}\label{lem:weak increasement}{\rm \cite[Theorem~10.12]{Zhang}} Let $a, b\in  \mathbb{R}^n$ be such that $a \prec_w b$. Then
$\sum\limits_{i=1}^nf(a_i)\leq \sum\limits_{i=1}^nf(b_i)$ for any increasing convex function $f$ on an interval containing all $a_i$ and $b_i$.
\end{lemma}

\begin{lemma}\label{lem:permutation is needed} {\rm \cite[Theorem~10.14]{Zhang}} Let $a, b\in  \mathbb{R}^n$ be such that $b$ is not a permutation of $a$. Then for any strictly increasing and strictly convex function $f$ that contains all the components of $a$ and $b$, we have
$$\sum\limits_{i=1}^nf(a_i)<\sum\limits_{i=1}^nf(b_i)\quad\mbox{whenever $a \prec_w b$.}$$
\end{lemma}

As a consequence of Lemmas~\ref{lem:weak increasement}--\ref{lem:permutation is needed} and Lemma~\ref{elemental observations}(ii), we have the following result.
\begin{lemma}\label{majorization of vector}  For every $T=(t_{ij})\in \mathbb{C}^{n\times n}$, we have
\begin{equation}\label{equ:trace bigger}
  \mbox{tr}\left[\left(I_n+TT^*\right)^{\frac{1}{2}}\right]\geq \sum_{i=1}^n\sqrt{1+|t_{ii}|^2}.
\end{equation}
Moreover, the equality above occurs only if $T$ is diagonal.
\end{lemma}
\begin{proof}Let $f(x)=\sqrt{1+x^2}$ for $x\in [0,+\infty)$.
It is clear that $f$ is a strictly increasing and strictly convex function on the interval $[0,+\infty)$, and
$$ \mbox{tr}\left[\left(I_n+TT^*\right)^{\frac{1}{2}}\right]= \sum_{i=1}^n f\big(\sigma_i(T)\big),$$
in which $\sigma_i(T)$ $(1\leq i\leq n)$ are the singular values of $T$.
Denote by
$$|d(T)|=\big(|t_{11}|, |t_{22}|,\cdots, |t_{nn}|\big),\quad \sigma(T)=\big(\sigma_1(T), \sigma_2(T),\cdots, \sigma_n(T) \big).$$
By \cite[Theorem~10.19]{Zhang} $|d(T)|\prec_w \sigma(T)$, which leads by Lemma~\ref{lem:weak increasement} to
 $$\sum_{i=1}^n f\big(\sigma_i(T)\big)\geq \sum_{i=1}^n f(|t_{ii}|).$$
This shows the validity of \eqref{equ:trace bigger}.

Suppose that \eqref{equ:trace bigger} turns out to be an equality.
Then by Lemma~\ref{lem:permutation is needed}
$|d(T)|=\sigma(T)W$ for some permutation matrix $W$, which implies that
$$\sum\limits_{i=1}^n |t_{ii}|^2=\sum\limits_{i=1}^n \sigma^2_i(T).$$
Hence, due to Lemma~\ref{elemental observations}(ii) $T$ is a diagonal matrix.
\end{proof}

Now, we turn to the description of the matched projection.
\begin{definition}\cite[Definition~3.1 and Theorem~3.2]{TXF} For each idempotent $Q\in \mathbb{C}^{n\times n}$, its matched projection is defined by
\begin{equation}\label{equ:final exp of P}m(Q)=\frac12\big(|Q^*|+Q^*\big)|Q^*|^\dag\big(|Q^*|+I_n\big)^{-1}\big(|Q^*|+Q\big),
\end{equation}
where $|Q^*|^\dag$ denotes the Moore-Penrose inverse of $|Q^*|$.
\end{definition}

\begin{remark}For every idempotent $Q$, it is known \cite[Theorem~3.2]{TXF} that $$|Q^*|^\dag=\left(P_{\mathcal{R}(Q)}P_{\mathcal{R}(Q^*)}P_{\mathcal{R}(Q)}\right)^\frac12.$$
\end{remark}
 Note that the matrix $m(Q)$ defined as above is a projection, which is equal to $Q$ whenever $Q$ is a projection. It is useful to derive block matrix representations for $m(Q)$. Given an idempotent $Q\in \mathbb{C}_r^{n\times n}$ with $1\le r<n$, let $U\in \mathbb{C}^{n\times n}$ be an arbitrary  unitary such that
\begin{equation}\label{decomposition of projection of Q}
P_{\mathcal{R}(Q)}=U^*\left(
                          \begin{array}{cc}
                            I_r & 0 \\
                            0 & 0 \\
                          \end{array}
                        \right)U.
\end{equation}
It follows that
\begin{equation}\label{decomposition of Q}Q=U^*\left(
                          \begin{array}{cc}
                            I_r & A \\
                            0 & 0 \\
                          \end{array}
                        \right)U\quad \mbox{
for some $A\in \mathbb{C}^{r\times (n-r)}$,}\end{equation}
which gives
$$|Q^*|=U^*\left(
             \begin{array}{cc}
               B & 0 \\
               0 & 0 \\
             \end{array}
           \right)U,$$
where \begin{equation}\label{equ:def of B}
B=(I_r+AA^*)^\frac12.
\end{equation}
Since $B\ge I_r$, we have  $\mbox{tr}(B)\ge \mbox{tr}(I_r)=r$. Therefore,
\begin{equation}\label{three traces are equal}\mbox{tr}(|Q^*|)\ge \mbox{rank}(Q)=\mbox{tr}(Q)=\mbox{tr}(Q^*)=\mbox{tr}\big(P_{\mathcal{R}(Q)}\big).\end{equation}
Furthermore, it can be concluded by Lemma~\ref{elemental observations}(i) that
\begin{equation}\label{equ:trace delete absolute value}\mbox{tr}(|Q^*|)=\mbox{tr}(Q)\Longleftrightarrow Q=P_{\mathcal{R}(Q)}\Longleftrightarrow \ \mbox{$Q$ is a projection}.\end{equation}

\begin{lemma}{\rm (cf.\,\cite[Theorem~3.1]{TXF})} Suppose that $Q$ is represented by \eqref{decomposition of Q}. Then
\begin{align}&\label{defn of widetilde Q wrt A B}
m(Q)=\frac12 U^*\left(\begin{array}{cc}
        (B+I_r)B^{-1} &  B^{-1}A\\
        A^*B^{-1} &  A^*\big[B(B+I_r)\big]^{-1}A \\
      \end{array}\right)U,
\end{align}
where $B$ is defined by \eqref{equ:def of B}.

\end{lemma}
\begin{proof}A simple use of \eqref{equ:final exp of P} and \eqref{decomposition of Q} yields the desired conclusion.
\end{proof}

We end this section by providing a formula for $\|m(Q)\|_F$.

\begin{lemma}\label{lem:f norm Q}For every idempotent $Q\in \mathbb{C}^{n\times n}$, we have
\begin{equation}\label{equ:f norm of mq and prq}
  \|m(Q)\|_F=\sqrt{\mbox{rank} (Q)}.
\end{equation}
\end{lemma}
\begin{proof}Let $$V=\frac{\sqrt{2}}{2}\big(|Q^*|+Q^*\big)\big(|Q^*|^\dag\big)^{\frac{1}{2}}\big(|Q^*|+I_n\big)^{-\frac{1}{2}}.$$
By \cite[Theorem~3.6]{TXF},  we have
$$m(Q)=VV^* \quad \mbox{and} \quad P_{\mathcal{R}(Q)}=V^*V.$$
Hence,
\begin{equation}\label{trace of m Q}\|m(Q)\|_F^2=\mbox{tr}\big(m(Q)\big)=\mbox{tr}\big(P_{\mathcal{R}(Q)}\big)=\mbox{rank} (Q).\end{equation}
This shows the validity of \eqref{equ:f norm of mq and prq}.
\end{proof}

\section{The minimum distance}\label{sec:positive answer}
In this section, we study the minimum distance from projections to an idempotent $Q$. It is helpful to begin with a formula for $\|m(Q)-Q\|_F$.
\begin{lemma}\label{lem:f norm m Q minus Q} For every idempotent $Q\in \mathbb{C}^{n\times n}$, we have
\begin{equation}\label{equ:F norm of mq-q}
  \|m(Q)-Q\|_F=\sqrt{\mbox{tr}(|Q^*|^2-|Q^*|)}.
\end{equation}
\end{lemma}
\begin{proof}  Since $Q$ is an idempotent and $\mathcal{R}(|Q^*|)=\mathcal{R}(Q)$,
we have
$$Q|Q^*|=|Q^*|\quad\mbox{and}\quad |Q^*|\cdot |Q^*|^\dag=|Q^*|^\dag \cdot |Q^*|=P_{\mathcal{R}(Q)}.$$
It follows that
$$Q(|Q^*|+Q^*)|Q^*|^\dag=|Q^*|(I_n+|Q^*|)|Q^*|^\dag=P_{\mathcal{R}(Q)}(I_n+|Q^*|),$$
which is combined with \eqref{equ:final exp of P} to get
\begin{equation*}
  Qm(Q)=\frac{1}{2}P_{\mathcal{R}(Q)}(|Q^*|+Q)=\frac{1}{2}(|Q^*|+Q).
\end{equation*}
Taking $*$-operation yields  $  m(Q)Q^*=\frac{1}{2}(|Q^*|+Q^*)$.
Therefore,
\begin{align*}
  \big( m(Q)-Q\big)\big(m(Q)-Q\big)^* & =m(Q)-Qm(Q)-m(Q)Q^*+QQ^* \\
   & =m(Q)-\frac{1}{2}(Q+Q^*)-|Q^*|+|Q^*|^2.
\end{align*}
In virtue of \eqref{three traces are equal} and \eqref{trace of m Q}, we have
$$\mbox{tr}\big[m(Q)-\frac{1}{2}(Q+Q^*)\big]=0.$$
Consequently,
$$\|m(Q)-Q\|_F^2=\mbox{tr}\big[\big( m(Q)-Q\big)\big(m(Q)-Q\big)^*\big]=\mbox{tr}(|Q^*|^2-|Q^*|).$$
So, the desired conclusion follows.
\end{proof}

Our first characterization of the minimum distance reads as follows.
\begin{theorem}\label{thm:existence}For every idempotent $Q\in \mathbb{C}^{n\times n}$, we have
\begin{equation}\label{m Q best}\|m(Q)-Q\|_F\le \|P-Q\|_F,\quad \forall\,P\in \mathcal{P}\left(\mathbb{C}^{n\times n}\right).\end{equation}
\end{theorem}
\begin{proof}Let $r=\mbox{rank}(Q)$. If $r=0$ or $r=n$, then $m(Q)=Q$, hence \eqref{m Q best} is trivially satisfied.

Suppose now that $1\le r<n$. Let $U\in \mathbb{C}^{n\times n}$ be a unitary such that
$Q$ is formulated by \eqref{decomposition of Q}. For any projection $P$ in $\mathbb{C}^{n\times n}$, from \cite[Theorem~2.14]{XY} we know that $P$ can be represented by
\begin{equation}\label{equ:2nd form of P}P=U^*\left(
                    \begin{array}{cc}
                      C & C^\frac12\big(I_r-C\big)^\frac12 U_0^* \\
                     U_0C^\frac12\big(I_r-C\big)^\frac12  & U_0\big(I_r-C\big)U_0^*+Q_0\\
                    \end{array}
                  \right)U,
\end{equation}
where  $C\in \mathbb{C}^{r\times r}$ is a positive contraction,
$U_0\in \mathbb{C}^{(n-r)\times r}$ is a partial isometry and $Q_0\in \mathbb{C}^{(n-r)\times (n-r)}$ is a projection satisfying
\begin{equation}\label{two pre-conditions on the decomposition}\mathcal{R}(U_0^*)=\mathcal{R}(C-C^2)\quad\mbox{and}\quad U_0^*Q_0=0.
\end{equation}
Let
\begin{equation}\label{defn of S}S=(Q+Q^*-I_n)P,\quad D=C^{\frac12}(I_r-C)^{\frac12}.\end{equation} Since $\mbox{tr}(PQ^*)=\mbox{tr}(Q^*P)$, we have
\begin{align*}\|P-Q\|_F^2=&\mbox{tr}\left[(P-Q)(P-Q)^*\right]=\mbox{tr}(QQ^*+P-QP)-\mbox{tr}(PQ^*)\\
=&\mbox{tr}(QQ^*-S)=\mbox{tr}(|Q^*|^2)-\mbox{tr}(S).
\end{align*}
Therefore, we may use \eqref{equ:F norm of mq-q} to get
\begin{align*}\|P-Q\|_F^2-\|m(Q)-Q\|_F^2=&\left[\mbox{tr}(|Q^*|^2)-\mbox{tr}(S)\right]-\left[\mbox{tr}(|Q^*|^2-|Q^*|)\right]\\
=&\mbox{tr}(|Q^*|)-\mbox{tr}(S)=\mbox{tr}\big[(I_r+AA^*)^\frac12\big]-\mbox{tr}(S).
\end{align*}
In view of \eqref{decomposition of Q}, \eqref{equ:2nd form of P} and \eqref{defn of S}, it is easy to verify that
$$S=U^*\left(
      \begin{array}{cc}
        C+AU_0 D & * \\
        * & A^*D U_0^*-U_0(I_r-C)U_0^*-Q_0 \\
      \end{array}
    \right)U,$$
which leads by $\mbox{tr}(Q_0)\ge 0$ to
\begin{align*}\mbox{tr}(S)=& \mbox{tr}\big[C+AU_0 D+A^*D U_0^*-U_0(I_r-C)U_0^*-Q_0\big]\\
\le&\mbox{tr}\big[C+AU_0 D+A^*D U_0^*-U_0(I_r-C)U_0^*\big].
\end{align*}
Utilizing $D=D^*$ and the observations
\begin{align*}&\mbox{tr}(A^*D U_0^*)=\overline{\mbox{tr}(U_0DA)}=\overline{\mbox{tr}(AU_0D)},\\
&\mbox{tr}\left[U_0(I_r-C)U_0^*\right]=\mbox{tr}\left[U_0^*U_0(I_r-C)\right],\end{align*}
we arrive at
\begin{align*}\mbox{tr}(S)\le& \mbox{tr}(C)-\mbox{tr}\left[U_0^*U_0(I_r-C)\right]+2\mbox{Re}\left[\mbox{tr}(AU_0D)\right].
\end{align*}
So, it is sufficient to prove that
\begin{equation}\label{remaining inequ}\mbox{tr}(C)-\mbox{tr}\left[U_0^*U_0(I_r-C)\right]+2\mbox{Re}\left[\mbox{tr}(AU_0D)\right]\le \mbox{tr}\left[(I_r+AA^*)^\frac12\right].
\end{equation}
For this, we choose a unitary $W\in \mathbb{C}^{r\times r}$ such that
\begin{equation}\label{equ:big C}C=W^*\cdot \mbox{diag}(c_1,c_2,\cdots,c_r)\cdot W,\end{equation}
in which $c_i\in [0,1]$ for all $i=1,2,\cdots,r$. Due to \eqref{defn of S}, we have
\begin{equation}\label{equ:big D}D=W^*\cdot \mbox{diag}(d_1,d_2,\cdots,d_r)\cdot W,\end{equation}
where
$d_i=\sqrt{c_i(1-c_i)}$  for $1\le i\le r$.
Taking a permutation if necessary, we may as well assume that
$d_1\ge d_2\ge \cdots\ge d_r\ge 0$.

\textbf{Case 1.}\quad  $d_r>0$. In this case, $\mbox{rank}(C-C^2)=\mbox{rank}(D^2)=r$, which implies that $U_0^*U_0=I_r$, since $U_0^*U_0$ is a projection satisfying the first equation in \eqref{two pre-conditions on the decomposition}. Consequently,
\eqref{remaining inequ} can be simplified as
\begin{equation}\label{reduced ineq}\mbox{tr}(2C-I_r)+2\mbox{Re}\left[\mbox{tr}(AU_0D)\right]\le \mbox{tr}\left[(I_r+AA^*)^\frac12\right].
\end{equation}
Let
$T=WAU_0W^*$ with $T=(t_{ij})_{1\le i,j\le r}$. Then
\begin{align*}\mbox{Re}\left[\mbox{tr}\big(AU_0D)\right]=&\mbox{Re}\left[\mbox{tr}(AU_0W^*\cdot \mbox{diag}(d_1,d_2,\cdots,d_r)\cdot W\big)\right]\\
=&\mbox{Re}\left[\mbox{tr}\big(T\cdot\mbox{diag}(d_1,d_2,\cdots,d_r)\big)\right]\\
=&\sum_{i=1}^r \mbox{Re}(t_{ii})\cdot d_i\le \sum_{i=1}^r |t_{ii}|\cdot d_i.
\end{align*}
Meanwhile,
\begin{align*}\mbox{tr}(2C-I_r)=&\mbox{tr}\left[W^*\big(2\cdot \mbox{diag}(c_1,c_2,\cdots,c_r)-I_r\big)W\right]\\
=&\mbox{tr}\big[2\cdot \mbox{diag}(c_1,c_2,\cdots,c_r)-I_r\big]=\sum_{i=1}^r(2c_i-1),
\end{align*}
and
\begin{align*}(2d_i)^2+(2c_i-1)^2=1\quad (1\le i\le r).
\end{align*}
As a result,  we may use Cauchy-Schwartz inequality and Lemma~\ref{majorization of vector} to obtain
\begin{align*}\mbox{tr}(2C-I_r)+2\mbox{Re}\left[\mbox{tr}\big(AU_0D)\right]\le & \sum_{i=1}^r \big(|t_{ii}|\cdot (2d_i)+(2c_i-1)\big)\\
\le & \sum_{i=1}^r \left(\sqrt{1+|t_{ii}|^2}\cdot \sqrt{(2d_i)^2+(2c_i-1)^2}\right)\\
= & \sum_{i=1}^r \sqrt{1+|t_{ii}|^2}\le \mbox{tr}\left[(I_r+TT^*)^\frac12\right].
\end{align*}
Observe that $U_0U_0^*\le I_{n-r}$, so we have
\begin{align*}TT^*=WAU_0W^*\cdot WU_0^*A^*W^*=WAU_0U_0^*A^*W^*\le WAA^*W^*,
\end{align*}
which in turn gives
$$I_r+TT^*\le W(I_r+AA^*)W^*,$$
hence by \cite[Proposition~1.3.8]{Pedersen} we have
$$(I_r+TT^*)^\frac12\le \left[W(I_r+AA^*)W^*\right]^\frac12=W(I_r+AA^*)^\frac12W^*.$$
It follows that
$$\mbox{tr}\left[(I_r+TT^*)^\frac12\right]\le \mbox{tr}\left[W(I_r+AA^*)^\frac12W^*\right]=\mbox{tr}\left[(I_r+AA^*)^\frac12\right].$$
This shows the validity of \eqref{reduced ineq}.

\textbf{Case 2.}\quad $d_s>0$ and $d_{s+1}=0$ for some $s$ with $1\le s\le r-1$. In this case,
$$D=W^*\cdot \mbox{diag}(d_1,d_2,\cdots,d_s,0_{r-s})\cdot W.$$
Note that the positivity of $D$ ensures $\mathcal{R}(D^2)=\mathcal{R}(D)$, so by \eqref{two pre-conditions on the decomposition} $U_0^*U_0$ is a projection whose range is equal to
that of $D$. Therefore,
$$U_0^*U_0=W^*\cdot \mbox{diag}(I_s,0_{r-s})\cdot W.$$
Let
\begin{equation*}
c_i^\prime=\left\{
    \begin{array}{ll}
      2c_i-1, & \hbox{$1\le i\le s$,} \\
      c_i, & \hbox{$s+1\le i\le r$,}
    \end{array}
  \right.
\quad
d_i^\prime=\left\{
    \begin{array}{ll}
      d_i, & \hbox{$1\le i\le s$,} \\
      0, & \hbox{$s+1\le i\le r$.}
    \end{array}
  \right.
\end{equation*}
Then $(c_i^\prime)^2+(2d_i^\prime)^2\le1$ for $1\le i \le r$, and
\begin{align*}\mbox{tr}(C)-\mbox{tr}\left[U_0^*U_0(I_r-C)\right]=\sum_{i=1}^r c_i-\sum_{i=1}^s (1-c_i)=\sum_{i=1}^r c_i^\prime.
\end{align*}
Let $T=(t_{ij})_{1\le i,j\le r}$ be defined as above. In virtue of
\begin{align*}
  AU_0D & =AU_0 W^*\cdot \mbox{diag}(d_1,d_2,\cdots , d_s,0_{r-s})\cdot W,
\end{align*}
we have
\begin{align*}
  \mbox{Re}\left[\mbox{tr}\big(AU_0D)\right] & =\mbox{Re}\left[\mbox{tr}\big(T\cdot \mbox{diag}(d_1,d_2,\cdots , d_s,0_{r-s})\big)\right] \\
   & =\sum_{i=1}^s\mbox{Re}(t_{ii})\cdot d_i\le \sum_{i=1}^s|t_{ii}|\cdot d_i= \sum_{i=1}^r|t_{ii}|\cdot d_i^\prime.
\end{align*}
Therefore, \eqref{remaining inequ} can be derived by using the same arguments employed in the derivation of \eqref{reduced ineq}.

\textbf{Case 3.}\quad $d_i=0$ for $1\le i\le r$. In this case $D=0$, so it can be concluded  by \eqref{defn of S} and \eqref{two pre-conditions on the decomposition} that $U_0^*U_0=0$ and $C$ is a projection. Hence,
\begin{equation}\label{c less less}\mbox{tr}(C)\le \mbox{tr}(I_r)\le \mbox{tr}\left[(I_r+AA^*)^\frac12\right].\end{equation}
Therefore, \eqref{remaining inequ} is also satisfied.  This completes the proof.
\end{proof}

Our second characterization of the minimum distance is stated as follows.
\begin{theorem}\label{thm:uniqueness} Let $Q\in \mathbb{C}^{n\times n}$ be an idempotent and $P\in\mathbb{C}^{n\times n}$ be a projection. If
$\|P-Q\|_F=\|m(Q)-Q\|_F$, then $P=m(Q)$.
\end{theorem}
\begin{proof} Suppose that $\|P-Q\|_F=\|m(Q)-Q\|_F$. It needs only to consider the case that $Q$ is not a projection. In this case, $Q\in \mathbb{C}_r^{n\times n}$ for some $1\le r<n$.  We follow the notations as that in the proof of Theorem~\ref{thm:existence}.

\textbf{Case 1.}\quad $d_r>0$. From Case 1 in the proof of Theorem~\ref{thm:existence}, it can be inferred  that
\begin{align*}&Q_0=0,\quad \mbox{Re}(t_{ii})\cdot d_i=|t_{ii}|\cdot d_i\quad (1\le i\le r),\\
&|t_{ii}|\cdot (2d_i)+(2c_i-1)=\sqrt{1+|t_{ii}|^2}\quad (1\le i\le r),\\
&\sum_{i=1}^r \sqrt{1+|t_{ii}|^2}=\mbox{tr}\left[(I_r+TT^*)^\frac12\right],\\
&TT^*=WAU_0U_0^*A^*W^*=WAA^*W^*\ \mbox{(by Lemma~\ref{elemental observations}(i))}.
\end{align*}
So $t_{ii}\ge 0$ for $1\le i\le r$\footnote{Actually, since $d_i$ is given by \eqref{equ:small c and d} and $d_i>0$, we have $t_{ii}>0$.},  and a simple application of Lemma~\ref{elemental observations}(iii) to the vectors $(t_{ii},1)$ and $\big(2d_i, 2c_i-1\big)$ in $\mathbb{R}^2$ yields
\begin{equation}\label{equ:small c and d} c_i=\frac{1+\sqrt{1+t_{ii}^2}}{2\sqrt{1+t_{ii}^2}},\quad d_i=\frac{t_{ii}}{2\sqrt{1+t_{ii}^2}}
\end{equation}for $1\le i\le r$.
Also, it can be concluded from Lemma~\ref{majorization of vector} that $T$ is a diagonal matrix.  So,  $T$ is in fact a semi-positive definite diagonal matrix.
Since $W$ is unitary and $WAU_0U_0^*A^*W^*=WAA^*W^*$, we have
\begin{equation}\label{expression for AA stra}AU_0U_0^*A^*=AA^*.\end{equation}
 Also, by definition we have  $T=WAU_0W^*$, so
$$AU_0=W^*TW=(W^*TW)^*=U_0^*A^*,$$
which is combined with \eqref{expression for AA stra} to get $AA^*=W^* T^2 W$, and thus
\begin{equation}\label{temp-001}|A^*|=W^*TW=AU_0.\end{equation}
Therefore, we may use the first equation above to obtain
$$W^*\cdot \mbox{diag}(t_{11}, t_{22}, \cdots, t_{rr})\cdot W=|A^*|.$$
This, together with \eqref{equ:big C}, \eqref{equ:big D} and \eqref{equ:small c and d}, yields
\begin{equation}\label{equ:big C+ and D+} C=\frac12 (B+I_r)B^{-1},\quad D=\frac12 B^{-1} |A^*|,\end{equation}
in which $B$ is defined by \eqref{equ:def of B}. Since $U_0U_0^*$ is a projection, we may rewrite \eqref{expression for AA stra}  as
\begin{equation*}\big[A(I_r-U_0U_0^*)\big]\big[A(I_r-U_0U_0^*)\big]^*=0,\end{equation*}
which happens only if $AU_0U_0^*=A$. It follows from  \eqref{temp-001} that
\begin{equation}\label{temp-002}|A^*| U_0^*=A.\end{equation}
Exploring \eqref{defn of S}, \eqref{equ:big C+ and D+} and \eqref{temp-002}, we arrive at
\begin{align*}&C^\frac12 (I_r-C)^\frac12 U_0^*=DU_0^*=\frac12 B^{-1} |A^*|U_0^*=\frac12 B^{-1} A,\\
&I_r-C=\frac12 (I_r-B^{-1})=\frac12 |A^*|\cdot \big[B(B+I_r)\big]^{-1}\cdot |A^*|\quad\mbox{(by \eqref{equ:def of B}) },\\
&U_0(I_r-C)U_0^*=\frac12 A^*\big[B(B+I_r)\big]^{-1}A.
\end{align*}
The arguments above  together with \eqref{defn of widetilde Q wrt A B} and \eqref{equ:2nd form of P} yield $P=m(Q)$.

\textbf{Case 2.}\quad $d_s>0$ and $d_{s+1}=0$ for some $s$ with $1\le s\le r-1$. By Case 2 in the proof of Theorem~\ref{thm:existence}, we have
\begin{align*}&Q_0=0,\quad \mbox{Re}(t_{ii})\cdot d_i=|t_{ii}|\cdot d_i\quad (1\le i\le s),\\
&|t_{ii}|\cdot (2d_i^\prime)+c_i^\prime=\sqrt{1+|t_{ii}|^2}\quad (1\le i\le r),\\
&\sum_{i=1}^r \sqrt{1+|t_{ii}|^2}=\mbox{tr}\left[(I_r+TT^*)^\frac12\right],\\
&TT^*=WAU_0U_0^*A^*W^*=WAA^*W^*.
\end{align*}
So $t_{ii}>0$ for $1\le i\le s$, \eqref{equ:small c and d} is  true for $1\le i\le s$, and
$$c_i=1,\quad t_{ii}=0\quad\mbox{($s+1\le i\le r$)}.$$
Therefore, \eqref{equ:small c and d} is valid for every $i$ with $1\le i\le r$. The rest of the proof is the same as that employed in Case 1.

\textbf{Case 3.}\quad $d_i=0$ for $1\le i\le r$. As is shown before, in this case we have $Q_0=0$, $D=U_0^*U_0=0$ and $C$ is a projection.
By assumption $\|P-Q\|_F=\|m(Q)-Q\|_F$, so actually \eqref{remaining inequ} turns out to be an equation. Therefore,
$$\mbox{tr}(C)=\mbox{tr}\left[(I_r+AA^*)^\frac12\right].$$
In view of \eqref{c less less}, the above equality is valid if and only if
$C=I_r$ and meanwhile $A=0$, which is exactly the case that $P=Q=P_{\mathcal{R}(Q)}$.
Hence, $P=m(Q)$ as desired.
\end{proof}

\section{The maximum distance and the intermediate value theorem}\label{sec:additional remark}
 The purpose of this section is to give positive answers to the remaining issues stated in Problem~\ref{problems}.
\begin{lemma}\label{lem:invariant}For every square matrix $A\in \mathbb{C}^{n\times n}$, the number $$\|P-A\|_F^2+\|I_n-P-A\|_F^2$$ is invariant with respect to the choice of the projection $P$ in $\mathbb{C}^{n\times n}$.
\end{lemma}
\begin{proof} For each projection $P\in\mathbb{C}^{n\times n}$, a simple computation yields
\begin{align*} &\|P-A\|_F^2=\mbox{tr}(P-PA^*-AP+AA^*),\\
&\|I_n-P-A\|_F^2=\mbox{tr}\left[(I_n-P)-(I_n-P)A^*-A(I_n-P)+AA^*\right].
\end{align*}
This shows that
$$\|P-A\|_F^2+\|I_n-P-A\|_F^2=\mbox{tr}(I_n-A^*-A+2AA^*),$$
which does not dependent on the choice of $P$.
\end{proof}

Our characterization of the maximum distance reads as follows.
\begin{theorem}\label{thm:maximum distance}For every idempotent $Q\in \mathbb{C}^{n\times n}$, we have
\begin{equation*}\label{m Q worst}\|P-Q\|_F\le \|I_n-m(Q)-Q\|_F,\quad \forall\,P\in \mathcal{P}\left(\mathbb{C}^{n\times n}\right).\end{equation*}
Moreover, the equality above occurs  only if $P=I_n-m(Q)$.
\end{theorem}
\begin{proof}For every projection $P\in \mathbb{C}^{n\times n}$, by Lemma~\ref{lem:invariant}, Theorems~\ref{thm:existence} and \ref{thm:uniqueness}, we have
\begin{align*}\|I_n-m(Q)-Q\|^2_F-\|P-Q\|^2_F=\|I_n-P-Q\|^2_F-\|m(Q)-Q\|_F^2\ge 0,
\end{align*}
and the number zero is obtained  only if $I_n-P=m(Q)$.
\end{proof}

To get a positive answer to Problem~\ref{problems}(iv), we need some additional lemmas.

\begin{lemma}\label{lem:two relationships with m q}{\rm \cite[Theorems~3.4 and 3.10]{TXF}} For every idempotent $Q\in \mathbb{C}^{n\times n}$, we have
\begin{equation*}\label{2 equs wrt m q}m(Q^*)=m(Q), \quad m(I_n-Q)=I_n-m(Q).\end{equation*}
\end{lemma}

\begin{lemma}\label{lem:homotopy-01}For every idempotent $Q\in \mathbb{C}^{n\times n}$, $P_{\mathcal{R}(Q)}$  and $m(Q)$ as projections are homotopy equivalent.
\end{lemma}
\begin{proof}It needs only to consider the case that $1\le r:=\mbox{rank}(Q)\le n-1$. Let $Q$ be represented by
\eqref{decomposition of Q}. For each $t\in [0,1]$, set
$$Q_t=U^*\left(
                          \begin{array}{cc}
                            I_r & tA \\
                            0 & 0 \\
                          \end{array}
                        \right)U.$$
Clearly, $Q_t(t\in [0,1])$ is a norm-continuous path of idempotents in $\mathbb{C}^{n\times n}$ which starts at $P_{\mathcal{R}(Q)}$ and ends at $Q$. This shows that  $P_{\mathcal{R}(Q)}$  and $Q$ as idempotents are homotopy equivalent. The same is also true for $Q$ and $m(Q)$ by  a direct use of \cite[Theorem~2.1]{TXF}, since the operator norm and the Frobenius norm on $\mathbb{C}^{n\times n}$ are equivalent. As  both  of $P_{\mathcal{R}(Q)}$ and $m(Q)$ are projections, there exists a norm-continuous path of projections from $P_{\mathcal{R}(Q)}$ to $m(Q)$ \cite[Proposition~4.6.3]{Blackadar}.
 \end{proof}

We are now in the position to prove the following intermediate value theorem.

\begin{theorem}\label{intermediate value theorem}For any idempotent $Q\in \mathbb{C}^{n\times n}$ with $Q\ne 0_n$ and
$Q\ne I_n$, we have $\mbox{ran}_Q=[\mbox{min}_Q, \mbox{max}_Q]$, where  $\mbox{ran}_Q$, $\mbox{min}_Q$ and $\mbox{max}_Q$
are defined by \eqref{defn of ran Q} and \eqref{Q min}, respectively.
\end{theorem}
\begin{proof}By Theorems~\ref{thm:existence} and \ref{thm:maximum distance}, we have
$$ \mbox{min}_Q=\|m(Q)-Q\|_F,\quad \mbox{max}_Q=\|I_n-m(Q)-Q\|_F.$$
For simplicity, we put
$$\lambda_1=\|P_{\mathcal{R}(Q)}-Q\|_F,\quad \lambda_2=\|I_n-P_{\mathcal{R}(Q)}-Q\|_F.$$
From \eqref{decomposition of projection of Q} and \eqref{decomposition of Q}, it is easily seen that
$$\lambda_1^2=\|A\|_F^2,\quad \lambda_2^2=\|A\|_F^2+n,$$
and thus $\lambda_2=\sqrt{\lambda_1^2+n}$. Therefore,
$[\mbox{min}_Q,\mbox{max}_Q]=J_1\cup J_2\cup J_3$,
where
\begin{equation*}J_1=\big[\|m(Q)-Q\|_F,\lambda_1\big],\quad J_2=[\lambda_1, \lambda_2],\quad J_3=\big[\lambda_2, \|I_n-m(Q)-Q\|_F\big].
\end{equation*}
So, it is reduced to show that $J_i\subseteq  \mbox{ran}_Q $ for $i=1,2,3$.

By Lemma~\ref{lem:homotopy-01}, there exists  a norm-continuous path of projections $P(t)(t\in [0,1])$ in $\mathbb{C}^{n\times n}$ such that $P(0)=m(Q)$ and $P(1)=P_{\mathcal{R}(Q)}$. Since the function: $t\to \|P(t)-Q\|_F $
is continuous on $[0,1]$, we see that $J_1\subseteq  \mbox{ran}_Q$.
In view of Lemma~\ref{lem:two relationships with m q}, we have
$$\|I_n-m(Q)-Q\|_F =\|m(I_n-Q^*)-Q\|_F.$$
Observe that
$$P_{\mathcal{R}(I_n-Q^*)}=P_{\mathcal{N}(Q^*)}=I_n-P_{\mathcal{R}(Q)},$$
so it can be concluded by Lemma~\ref{lem:homotopy-01} that $m(I_n-Q^*)$ and  $I_n-P_{\mathcal{R}(Q)}$ are homotopy equivalent. Hence, similar reasoning shows that $J_3\subseteq  \mbox{ran}_Q$.

Now, we turn to prove that $J_2\subseteq  \mbox{ran}_Q$. By assumption we have
$1\le r:=\mbox{rank}(Q)\le n-1$. In virtue of
$$\|(I_n-P)-(I_n-Q^*)\|_F=\left\|\left[(I_n-P)-(I_n-Q^*)\right]^*\right\|_F=\|P-Q\|_F$$
for every projection $P\in \mathbb{C}^{n\times n}$, and
\begin{align*}&\|P_{\mathcal{R}(I_n-Q^*)}-(I_n-Q^*)\|_F=\left\|\left[P_{\mathcal{R}(I_n-Q^*)}-(I_n-Q^*)\right]^*\right\|_F=\lambda_1,\\
&\mbox{rank}(I_n-Q^*)=\mbox{rank}(I_n-Q)=n-r,
\end{align*}
we may as well assume that $r\le n-r$. Let $P_{\mathcal{R}(Q)}$ and $Q$ be represented by \eqref{decomposition of projection of Q} and \eqref{decomposition of Q}, respectively. For each $t\in [0,1]$, let
\begin{align*}&c_t=\cos\left(\frac{\pi t}{2}\right), \quad s_t=\sin\left(\frac{\pi t}{2}\right), \quad P_t=\left(
               \begin{array}{cc}
                 c^2_t I_r & s_tc_t I_r   \\
                 s_tc_t I_r & s_t^2 I_r \\
                  \end{array}
             \right)\in \mathbb{C}^{2r\times 2r},\\
&P_t^{(0)}=\left(
               \begin{array}{c|c}
                 P_t &  \\\hline
                 &  0_{n-2r} \\
               \end{array}
             \right),\quad P_t^{(1)}=\left(
               \begin{array}{c|cc}
                 P_t & & \\\hline
                  &  1 & \\
                   & & 0_{n-2r-1}
               \end{array}
             \right),\\
&P_t^{(2)}=\left(
               \begin{array}{c|cc}
                 P_t & & \\\hline
                  &  I_2 & \\
                   & & 0_{n-2r-2}
               \end{array}
             \right),\cdots, P_t^{(n-2r)}=\left(
               \begin{array}{c|c}
                 P_t &  \\\hline
                  &  I_{n-2r} \\
                   \end{array}
             \right).
\end{align*}
For each $i\in\{0,1,2,\cdots, n-2r\}$, it is clear that $P_t^{(i)}$ $(t\in [0,1])$ is a norm-continuous path of projections. Hence, as is observed before,  a continuous function can be induced as
$$f_i(t)=\big\|P_t^{(i)}-Q\big\|_F\quad  (t\in [0,1]).$$
A simple calculation shows that
\begin{align*}&f_0(0)=\|A\|_F=\lambda_1,\quad f_0(1)=\sqrt{\lambda_1^2+2r},\\
&f_1(0)=\sqrt{\lambda_1^2+1},\quad f_1(1)=\sqrt{\lambda_1^2+2r+1},\\
&f_2(0)=\sqrt{\lambda_1^2+2},\quad f_2(1)=\sqrt{\lambda_1^2+2r+2},\\
& \qquad\quad \vdots\\
&f_{n-2r}(0)=\sqrt{\lambda_1^2+n-2r},\quad f_{n-2r}(1)=\sqrt{\lambda_1^2+2r+(n-2r)}=\lambda_2.
\end{align*}
Since $2r>1$, we see that
 $$[\lambda_1,\lambda_2]=\bigcup_{i=0}^{n-2r}\big[f_i(0),f_i(1)\big]\subseteq \mbox{ran}_Q.\qedhere$$
\end{proof}

\section{Some Frobenius  norm bounds}
In this section, we study the Frobenius  norm bounds of $P-Q$ and $P-Q^\dag$ under the restriction of $PQ=Q$ on  a projection $P$ and an idempotent $Q$.

\begin{theorem}\label{thm:app lower-01} Let $P\in \mathbb{C}^{n\times n}$ be a projection and $Q\in \mathbb{C}^{n\times n}$ be an idempotent such that $PQ=Q$. Then
\begin{equation}\label{equ:f norm inequ with p q mq}
  \|P-Q\|_{F}\ge \lambda_{P,Q},
\end{equation}
where
\begin{equation}\label{equ:bounded}
\lambda_{P,Q}=\sqrt{\|m(Q)-Q\|_F^2+\mbox{rank}(P)-\mbox{rank}(Q)}.
\end{equation}
Moreover, the equality above occurs if and only if $Q$ is a projection.
\end{theorem}
\begin{proof}
Denote by $\mbox{rank}(Q)=r$. If $r=0$, then $Q=0$, so \eqref{equ:f norm inequ with p q mq} is trivially satisfied.
If $r=n$, then $Q=I_n$, which yields $P=I_n$, since by assumption $PQ=Q$. Therefore, \eqref{equ:f norm inequ with p q mq} is also trivially satisfied.

Now, we suppose that $1\le r< n$ and $PQ=Q$. In this case, we have $\mbox{rank}(P)\ge \mbox{rank}(Q)$. Rewrite $QQ^*$ as
$QQ^*=PQQ^*P$, which implies that
\begin{equation*}\label{equ:range q include range p}
P|Q^*|=|Q^*|=|Q^*|P.
\end{equation*}
As a result,
\begin{equation}\label{equ:pre-squared}QQ^*-2|Q^*|+P=(|Q^*|-P)^2.
\end{equation}
Meanwhile, since
$\mbox{tr}(QP)=\mbox{tr}(PQ)$ and $Q^*P=(PQ)^*=Q^*$, it can be derived from  \eqref{three traces are equal} that
$$\mbox{tr}(QP)=\mbox{tr}(Q)=\mbox{tr}(Q^*P).$$
Hence,
\begin{equation*}\mbox{tr}(S)=2\cdot \mbox{tr}(Q)-\mbox{tr}(P),\end{equation*} where $S$ is defined by \eqref{defn of S}.
As is shown in the proof of Theorem~\ref{thm:existence}, we have
$$\|P-Q\|_F^2=\mbox{tr}(QQ^*-S),\quad \|m(Q)-Q\|_F^2=\mbox{tr}(QQ^*-|Q^*|).$$
It follows that
\begin{align}\label{squared F norm-01} &\|P-Q\|_F^2=\mbox{tr}(QQ^*)-2\cdot \mbox{tr}(Q)+\mbox{tr}(P),\\
\label{squared lam-01}&\lambda_{P,Q}^2=\mbox{tr}(QQ^*)-\mbox{tr}(|Q^*|)+\mbox{tr}(P)-\mbox{tr}(Q).
\end{align}
So, we may use \eqref{three traces are equal} and \eqref{equ:trace delete absolute value} to conclude that
$\|P-Q\|_F^2\ge \lambda_{P,Q}^2$ and $$\|P-Q\|_F^2=\lambda_{P,Q}^2\Longleftrightarrow Q=P_{\mathcal{R}(Q)}.$$
So, the desired conclusion follows.
\end{proof}

\begin{theorem}\label{thm:app up-01}Under the condition of  Theorem~\ref{thm:app lower-01}, we have
\begin{equation}\label{equ:f norm inequ with p q mq++}
  \|P-Q\|_{F}\le \sqrt{2}\cdot \lambda_{P,Q},
\end{equation}
where $\lambda_{P,Q}$ is defined by \eqref{equ:bounded}. Moreover, the equality above occurs if and only if $P=Q=P_{\mathcal{R}(Q)}$.
\end{theorem}
\begin{proof} To avoid the triviality, we only consider the case that $1\le \mbox{rank}(Q)<n$ and $PQ=Q$. Since
$|Q^*|$ and $P$ are both Hermitian, we may use \eqref{squared lam-01}, \eqref{squared F norm-01}  and \eqref{equ:pre-squared} to obtain
$$2\cdot \lambda_{P,Q}^2-\|P-Q\|_F^2=\mbox{tr}\left[(|Q^*|-P)^2\right]=\big\| |Q^*|-P\big\|_F^2.$$
The above equalities together with \eqref{three traces are equal} and \eqref{equ:trace delete absolute value} yield the conclusion.
\end{proof}

As an application of Theorem~\ref{thm:app up-01}, an upper of $ \|P-Q^\dag \|_{F}$ can be derived as follows.
\begin{theorem}\label{thm:app up-01 point 5}Under the condition of  Theorem~\ref{thm:app lower-01}, we have
\begin{equation}\label{equ:f norm inequ with p q mq+++}
  \|P-Q^\dag \|_{F}\le \sqrt{2}\cdot \sqrt{1+\|Q^\dag Q-QQ^\dag\|_2^2}\cdot \lambda_{P,Q},
\end{equation}
where $\lambda_{P,Q}$ is defined by \eqref{equ:bounded}. Moreover, the equality above occurs if and only if $P=Q=P_{\mathcal{R}(Q)}$.
\end{theorem}
\begin{proof}Since $Q$ is an idempotent, we have
$$Q^\dag=Q^\dag Q\cdot QQ^\dag=P_{\mathcal{R}(Q^*)}P_{\mathcal{R}(Q)},$$
which gives $\|Q^\dag\|_2\le 1$. A direct use of \cite[Theorem~1.2]{Koliha} yields
\begin{align*}P-Q^\dag=&\Omega_1+\Omega_2+\Omega_3,
\end{align*}
where
\begin{align*}&\Omega_1=-P(P-Q)Q^\dag, \quad \Omega_2=(I_n-P)(P-Q)^*(Q^\dag)^*Q^\dag,\\
&\Omega_3=P(P-Q)^*(I_n-QQ^\dag).
\end{align*}
With the decomposition of $P-Q^\dag$ as above, it is clear that
$$\|P-Q^\dag\|_F^2=\sum_{i=1}^3
\|\Omega_i\|_F^2.$$
Since $Q^2=Q$ and by assumption $PQ=Q$, we have
$(P-Q)QQ^\dag=0$. Therefore,
\begin{align*}\|\Omega_1\|_F=&\|P(P-Q)Q^\dag Q\cdot Q^\dag\|_F\le \|P\|_2\cdot \|(P-Q)Q^\dag Q\|_F\cdot \|Q^\dag\|_2\\
\le &\|(P-Q)Q^\dag Q\|_F=\|(P-Q)(Q^\dag Q-QQ^\dag)\|_F\\
\le & \|P-Q\|_F\cdot \|Q^\dag Q-QQ^\dag\|_2.
\end{align*}
Meanwhile,
\begin{align*}\|\Omega_2\|_F\le & \|(I_n-P)(P-Q)^*\|_F\cdot \|Q^\dag\|_2^2\le \|(I_n-P)(P-Q)^*\|_F\\
=&\big\|\left[(P-Q)(I_n-P)\right]^*\big\|_F=\|(P-Q)(I_n-P)\|_F,
\end{align*}
and
\begin{align*}\|\Omega_3\|_F=&\|\Omega_3^*\|_F=\|(I_n-QQ^\dag)(P-Q)P\|_F\\
\le &\|I_n-QQ^\dag\|_2\cdot  \|(P-Q)P\|_F\le \|(P-Q)P\|_F.\end{align*}
Consequently,
$$\|\Omega_2\|_F^2+\|\Omega_3\|_F^2\le \|P-Q\|_F^2.$$
The above arguments together with  \eqref{equ:f norm inequ with p q mq++} yield
\begin{align*}\|P-Q^\dag\|_F\le & \sqrt{1+\|Q^\dag Q-QQ^\dag\|_2^2}\cdot \|P-Q\|_F\\
\le &\sqrt{2}\cdot \sqrt{1+\|Q^\dag Q-QQ^\dag\|_2^2}\cdot \lambda_{P,Q}.
\end{align*}
Moreover, from the above derivation we see that if \eqref{equ:f norm inequ with p q mq+++} becomes an equality, then so does for
\eqref{equ:f norm inequ with p q mq++}. Hence,  $P=Q=P_{\mathcal{R}(Q)}$ as is desired.
\end{proof}

\begin{remark}Since both of $Q^\dag Q$ and $QQ^\dag$ are projections, the well-known Krein-Krasnoselskii-Milman equality \cite{Kato} indicates that $\|Q^\dag Q-QQ^\dag\|_2\le 1$.
\end{remark}

For each idempotent $Q$, it is clear that $P_{\mathcal{R}(Q)}Q=Q$. As is observed in  \eqref{three traces are equal},  we have $\mbox{rank}(P)=\mbox{rank}(Q)$ in the case that
$P=P_{\mathcal{R}(Q)}$. So, in view of Theorems~\ref{thm:app lower-01} and \ref{thm:app up-01}, a bilateral inequality can be derived immediately as follows.

\begin{theorem}\label{thm:the range projection} For every idempotent $Q\in \mathbb{C}^{n\times n}$, we have
\begin{equation}\label{equ:optimal coefficient}
  \frac{\sqrt{2}}{2}\cdot \|P_{\mathcal{R}(Q)}-Q\|_F  \le \|m(Q)-Q\|_F\le \|P_{\mathcal{R}(Q)}-Q\|_F.
\end{equation}
Moreover, each equality above occurs if and only if $Q$ is a projection.
\end{theorem}

\begin{remark}The second inequality in \eqref{equ:optimal coefficient} can also be obtained by a direct use of Theorem~\ref{thm:existence}.
\end{remark}

\begin{remark}\label{rem:optimal numbers} By Theorem~\ref{thm:the range projection}, there exist constants $\alpha_1, \alpha_2\in (0,+\infty)$ such that
\begin{equation}\label{equ:calimed sharpness max min}\alpha_1\cdot \left\|P_{\mathcal{R}(Q)}-Q\right\|_F\le \left\|m(Q)-Q\right\|_F\le \alpha_2 \cdot \left\|P_{\mathcal{R}(Q)}-Q\right\|_F
\end{equation}
for any $n\in \mathbb{N}$ and every idempotent $Q\in\mathbb{C}^{n\times n}$.
For each $a>0$, let $Q_a$ be the idempotent  given by
\begin{equation*}\label{defn of Q a}Q_a=\left(
        \begin{array}{cc}
          1 & a \\
          0 & 0 \\
        \end{array}
      \right)\in \mathbb{C}^{2\times 2}.\end{equation*}
Evidently, $\|P_{\mathcal{R}(Q_a)}-Q_a\|_F=a$ and by \eqref{equ:F norm of mq-q} we have
$$\left\|m(Q_a)-Q_a\right\|_F=\sqrt{1+a^2-\sqrt{1+a^2}}.$$
Therefore,
\begin{equation*}\label{two limits}\lim_{a\to 0^+}\frac{\left\|m(Q_a)-Q_a\right\|_F}{\left\|P_{\mathcal{R}(Q_a)}-Q_a\right\|_F}=\frac{\sqrt{2}}{2}\quad\mbox{and}\quad \lim_{a\to +\infty}\frac{\left\|m(Q_a)-Q_a\right\|_F}{\left\|P_{\mathcal{R}(Q_a)}-Q_a\right\|_F}=1.\end{equation*}
The limitations above together with \eqref{equ:optimal coefficient} indicate that $\frac{\sqrt{2}}{2}$ and $1$ are the optimal numbers of $\alpha_1$ and $\alpha_2$ that ensure the validity of \eqref{equ:calimed sharpness max min} for every idempotent $Q$.

\end{remark}

Based on Theorems~Theorem~\ref{thm:existence} and \ref{thm:the range projection}, an additional upper bound of
 $\|P_{\mathcal{R}(Q)}-Q\|_F$
can be put forward as follows.

\begin{theorem}\label{thm:ist app}For every idempotent $Q\in \mathbb{C}^{n\times n}$ and every projection $P\in \mathbb{C}^{n\times n}$, we have
\begin{equation*}
  \|P_{\mathcal{R}(Q)}-Q\|_F\le  \sqrt{2}\cdot \|P-Q\|_F.
\end{equation*}
Moreover, the equality above occurs if and only if $P=Q$.
\end{theorem}
\begin{proof}We may combine Theorems~\ref{thm:the range projection} and \ref{thm:existence} to conclude that
\begin{equation*}\|P_{\mathcal{R}(Q)}-Q\|_F\le \sqrt{2}\cdot \|m(Q)-Q\|_F\le  \sqrt{2}\cdot \|P-Q\|_F.
\end{equation*}
So, if $\|P_{\mathcal{R}(Q)}-Q\|_F=\sqrt{2}\cdot \|P-Q\|_F$, then $\|P_{\mathcal{R}(Q)}-Q\|_F=\sqrt{2}\cdot \|m(Q)-Q\|_F$,
hence by Theorem~\ref{thm:the range projection} $Q=P_{\mathcal{R}(Q)}$, and thus $\sqrt{2}\cdot \|P-Q\|_F=0$, which happens only if $P=Q$.
\end{proof}

\begin{remark}It is obvious that $PQ=Q$ for $P=AA^\dag$ and $Q=AX$, where $X$ denotes an arbitrary generalized inverse of a  matrix $A$, which includes  the weighted Moore-Penrose inverse \cite{QXZ}, the Drazin inverse \cite{XSW}  and so on. So, results  obtained as above can be applied to deal with the Frobenius norm estimations for $AA^\dag-AX$.
\end{remark}

Taking the Drazin inverse as an example, in what follows we derive two  Frobenius norm bounds for $ AA^\dag-AA^d$, where $A^d$ denotes the Drazin inverse of a square matrix $A$. It is well-known that
$$\mathcal{R}(AA^d)=\mathcal{R}(A^d)=\mathcal{R}(A^k),$$ where $k=\mbox{ind}(A)$, which is called the Drazin index of $A$.
It may happen that $AA^d$ is a projection, whereas $A^d\ne A^\dag$. For instance, if  $A$ is a non-zero nilpotent matrix, then $A^d=0$, while $A^\dag\ne 0$.
However, if $AA^d=AA^\dag$, then it  can be shown easily that $A^d$ and $A^\dag$ must be the same.

We may apply Theorems~\ref{thm:app lower-01} and \ref{thm:app up-01} to get the following norm bounds.

\begin{corollary}\label{thm:drazin lower bound} For every square matrix $A\in \mathbb{C}^{n\times n}$, we have
\begin{equation*}\label{equ:drazin low-01}\|AA^\dag-AA^d\|_F\ge \sqrt{\|m(AA^d)-AA^d\|_F^2+\mbox{rank}(A)-\mbox{rank}(A^d)}.
\end{equation*}
Moreover, the equality above occurs if and only if $AA^d$ is a projection.
\end{corollary}

\begin{corollary}\label{thm:drazin upper bound} For every square matrix $A\in \mathbb{C}^{n\times n}$, we have
\begin{equation}\label{equ:drazin upp-012}\|AA^\dag-AA^d\|_F\le \sqrt{2}\cdot \sqrt{\|m(AA^d)-AA^d\|_F^2+\mbox{rank}(A)-\mbox{rank}(A^d)}.
\end{equation}
Moreover, the equality above occurs if and only if $A^d=A^\dag$.
\end{corollary}

 We end this paper by showing the optimality of $\sqrt{2}$  associated with \eqref{equ:drazin upp-012} in the case that $\mbox{rank}(A)>\mbox{rank}(A^d)$. For this, we put
$$A_{n,a}=\left(
            \begin{array}{cc}
              I_{2n} & aI_{2n} \\
              0_{2n} & J_{n} \\
            \end{array}
          \right), \quad Q_{n,a}=\left(
      \begin{array}{cc}
        I_{2n} & a(I_{2n}+J_{n}) \\
        0_{2n} & 0_{2n} \\
      \end{array}
    \right)
$$
for each $n\in\mathbb{N}$ and $a\in (0,+\infty)$, where $J_{n}$ is a nilpotent defined by
$$J_{n}=\left(
    \begin{array}{c|c}
              \begin{array}{cc}
          0 & 1 \\
          0 & 0 \\
        \end{array}
            &  \\\hline
       & 0_{2n-2} \\
    \end{array}
  \right)\in \mathbb{C}^{2n\times 2n}.
$$
It is routine to verify that
$$A_{n,a}Q_{n,a}=Q_{n,a}A_{n,a}=Q_{n,a}=Q_{n,a}^2=A_{n,a}^2,$$
which means that
$Q_{n,a}=A_{n,a}^d$ such that $\mbox{ind}(A_{n,a})=2$. Hence,
$$A_{n,d}A^d_{n,a}=A_{n,a}Q_{n,a}=Q_{n,a}.$$
Since $\mathcal{R}(A_{n,a})=\mathbb{C}^{2n}\oplus \mathcal{R}(J_n)$, we have
$$A_{n,a}A_{n,a}^\dag=\left(
                        \begin{array}{cc}
                          I_{2n} &  \\
                           & J_nJ_n^\dag \\
                        \end{array}
                      \right)=\left(
                                \begin{array}{cc}
                                  I_{2n+1} &  \\
                                   & 0_{2n-1} \\
                                \end{array}
                              \right).
$$
Consequently, $\mbox{rank}(A_{n,a})-\mbox{rank}(A_{n,a}^d)=1$ and
$$\|A_{n,a}A_{n,a}^\dag-A_{n,a}A_{n,a}^d\|_F^2=\|A_{n,a}A_{n,a}^\dag-Q_{n,a}\|_F^2=1+(2n+1)a^2.$$
An easy calculation yields
$$Q_{n,a}Q_{n,a}^*=\left(
                     \begin{array}{c|c}
                                                \begin{array}{cc}
                           1+2a^2 & a^2 \\
                           a^2 & 1+a^2 \\
                         \end{array}
                                               &  \\\hline
                        &
                             \begin{array}{cc}
                               (1+a^2)I_{2n-2} & \\
                                & 0_{2n} \\
                             \end{array}

                        \\
                     \end{array}
                   \right),$$
so the non-zero eigenvalues  of $Q_{n,a}$  read as
$$\lambda_{1,2}=\frac{(2+3a^2)\pm\sqrt{5}a^2}{2},\quad  \lambda_i=1+a^2 \quad (3\le i\le 2n).$$
Utilizing \eqref{equ:F norm of mq-q} gives
\begin{align*}\|m(Q_{n,a})-Q_{n,a}\|_F^2=&\sum_{i=1}^{2n}(\lambda_i-\sqrt{\lambda_i})\\
=&b+(2n-2)\left[(1+a^2)-\sqrt{1+a^2}\right],
\end{align*}
in which
$$b=\lambda_1-\sqrt{\lambda_1}+\lambda_2-\sqrt{\lambda_2}.$$
Suppose now that $\alpha_3$ is a constant such that
\begin{equation*}\|A_{n,a}A_{n,a}^\dag-Q_{n,a}\|_F\le \alpha_3\cdot \sqrt{\|m(Q_{n,a})-Q_{n,a}\|_F^2+1}
\end{equation*}
for all  $n\in \mathbb{N}$ and $a\in (0,+\infty)$. Then
\begin{align*}\alpha_3^2\ge &\lim_{a\to 0^+}\lim_{n\to \infty}\frac{\|A_{n,a}A_{n,a}^\dag-Q_{n,a}\|_F^2}{\|m(Q_{n,a})-Q_{n,a}\|_F^2}\\
=&\lim_{a\to 0^+}\lim_{n\to \infty}\frac{1+(2n+1)a^2}{b+(2n-2)\left[(1+a^2)-\sqrt{1+a^2}\right]}\\
=&\lim_{a\to 0^+}\frac{a^2}{(1+a^2)-\sqrt{1+a^2}}=2.
\end{align*}
This shows the optimality of $\sqrt{2}$.

\vspace{2ex}

\noindent\textbf{Acknowledgement}

\vspace{2ex}

The authors thank Professor Fuzhen Zhang  for bring our attention to the theory of the weak majorization, and for his contribution to the proofs  of Lemmas~\ref{elemental observations} and \ref{majorization of vector}.




\end{document}